\documentclass {amsart}
\usepackage{graphicx}
\usepackage{amssymb}
\usepackage{amsmath}
\newtheorem{theorem}{Theorem}[section]

\newtheorem{remark}[theorem]{Remark}
\newtheorem{definition}{Definition}[section]

\newenvironment{pot2.1}{{\bf Proof.}}{\hfill\fbox{}\par\vspace{.2cm}}

\numberwithin{equation}{section}

\def\charf {\mbox{{\text 1}\kern-.24em {\text l}}}

\def\bea{\begin{eqnarray*}}
\def\eea{\end{eqnarray*}}
\def\be{\begin{eqnarray}}
\def\ee{\end{eqnarray}}

\begin{document}

\title[Myers-type theorems]{A note on generalized Myers-type theorems for h-almost Ricci tensors and generalized quasi-Einstein tensors}
%    Remove any unused author tags.

%    author one information
\author{Sanghun Lee}
\address{Department of Mathematics, Chung-Ang University, 84 Heukseok-ro, Dongjak-gu, Seoul, Republic of Korea}
\curraddr{}
\email{kazauye@cau.ac.kr}

%    author two information

% \thanks{This research was supported by Basic Science Research Program through the National Research Foundation of Korea (NRF) funded by the Ministry of Education (2019R1A6A3A13096187)}



\date{}

\dedicatory{}

\begin{abstract}
 In this paper, we prove some compactness theorems of Myers, Ambrose, and Galloway for complete Riemannian manifold in the concept of $h$–almost Ricci tensors and generalized quasi–Einstein tensors. Also, we extend the previous theorems when $h$ has at most linear growth in the distance function.
\end{abstract}

\maketitle
\section{Introduction}
 The concept of $h$-almost Ricci solitons is intoduced by Gomes, Wang, and Xia \cite{GWX}. This soliton is a natural extension of an almost Ricci soliton \cite{BR, PR}. An almost Ricci soliton is an $n$-dimensional Riemannian manifold $(M,g)$ with a vector field $V$ on $M$ and a soliton function $\lambda: M \rightarrow \mathbb{R}$ satisfying
\be \label{eq1.1}
Ric_{\, V} = \lambda g.
\ee
Here,
\bea
Ric_{\, V}:= Ric + \frac{1}{2}\mathcal{L}_{V}g,
\eea
where $\mathcal{L}$ denotes the Lie derivative. $Ric_{\, V}$ is called Bakry–Emery Ricci tensor and it is related to diffusion processess (see \cite{BE}). The Bakry–Emery Ricci tensor is studied in \cite{BE} and occurs naturally in many different subjects (cf. \cite{FG, ML, LO, GW}).  Recently, Gomes, Wang, and Xia come up with the following definition, which is a generalization of the almost Ricci soliton.
\begin{definition} [\cite{GWX}]
\emph{An} $h$-almost Ricci soliton \emph{is an $n$-dimensional Riemannian manifold $(M,g)$ with a vector field $V$ on $M$ and a soliton function $\lambda: M \rightarrow \mathbb{R}$ and a function $h: M \rightarrow \mathbb{R}$ which are smooth and satisfy the equation:}
\be \label{eq1.2}
Ric^{\, h}_{\, V} = \lambda g.
\ee
{\rm Here}, \\
\be \label{eq1.3}
Ric^{\, h}_{\, V}:= Ric + \frac{h}{2}\mathcal{L}_{V}g.
\ee
\end{definition}
We say that $Ric^{\, h}_{\, V}$ is an $\textit{h–almost Ricci tensor}$. When $V =\nabla u$ for some smooth function $u: M \rightarrow \mathbb{R}$, we call this a gradient $h$–almost Ricci soliton with a potential function $u$. In this case, the equation (\ref{eq1.2}) can be written as
\be \label{eq1.4}
Ric + h \, {\rm Hess}\, u = \lambda g,
\ee
where ${\rm Hess}$ $u$ denotes the Hessian of $u$. 

 It should be mentioned that an $h$-almost Ricci soliton is $\textit{expanding}$, $\textit{steady}$, or $\textit{shrinking}$ if $\lambda$ is negative, zero, or positive, respectively, on $M$. When $\mathcal{L}_{V}g = cg$ for some constant $c$, an $h$-almost Ricci soliton is said to be trivial. Otherwise it is nontrivial. We remark that the traditional Ricci soliton is a $1$-Ricci soliton with constant $\lambda$. Moreover, $1$-almost Ricci soliton is just the almost Ricci soliton. It is said that $h$ has defined to be signal if either $h>0$ on $M$ or $h<0$ on $M$.

 In \cite{MA} Maschler studied equation (\ref{eq1.4}) and he referred to equation (\ref{eq1.4}) as Ricci-Hessian equation. The Ricci-Hessian equation is related to a new class of Riemannian metrics, introduced by Catino \cite{CA}, which are natural generalizations of Einstein metrics. In more detail, he called a generalized quasi-Einstein manifold, if there are smooth functions $f$, $\lambda$, $\mu$ on $M$ satisfying
\be \label{eq1.5}
Ric_{f}^{\, \mu} = \lambda g.
\ee
Here,
\be \label{eq1.6}
Ric_{f}^{\, \mu} := Ric + \mbox{\rm Hess}\, f - \mu \, df \otimes df
\ee
and we say that $Ric_{f}^{\, \mu}$ is a $\textit{generalized quasi–Einstein tensor}$. When $\mu = \frac{1}{m}$, where $m$ is a positive integer, the above generalized quasi–Einstein manifold is called a generalized $m$-quasi–Einstein manifold (cf. \cite{BA}) and simply $m$-quasi–Einstein manifold when $\lambda$ is constant. It has been proved in \cite{JS, KK} that $m$-quasi–Einstein manifolds are directly related to the warped product Einstein manifolds. 

 The purpose of this paper is to investigate which geometric and topological results for manifolds with a lower bound on the Ricci tensor extend to smooth metric measure spaces with $h$–almost Ricci tensor or generalized quasi–Einstein tensor bounded below. Myers studied this problem in \cite{MY}. More precisely, Myers stated that if an $n$-dimensional complete Riemannian manifold $(M,g)$ satisfies $Ric \geq (n-1)H$ with $H>0$, then $M$ is compact and $\mbox{\rm diam(M)} \leq \frac{\pi}{\sqrt{H}}$. Moreover, the fundamental group $\pi_{1}(M)$ is finite. This theorem has been widely generalized in various directions by many authors \cite{AB, CGM, FG, GL, ML, CS, HT, LF, GW}. In \cite{AB}, Ambrose first generalized. Ambrose replaced lower bound on the Ricci curvature with an integral condition on the Ricci curvature.

\begin{theorem} [\cite{AB}] \label{thm1.1}
 Let $M$ be an $n$-dimensional complete Riemannian manifold. Suppose that there exists some point $p \in M$ for which every geodesic $\gamma: [0,\infty) \rightarrow M$ emanating from $p$ satisfies
\[ \int^{\infty}_{0} Ric(\gamma'(s),\gamma'(s)) ds = \infty. \]
Then $M$ is compact.
\end{theorem}

 On the other hand, motivated by relativistic cosmology, Galloway \cite{GL} proved compactness theorem by perturbing the positive lower bound on the Ricci curvature by the derivative in the radial direction of some bounded function.

\begin{theorem} [\cite{GL}] \label{thm1.2}
 Let $M$ be an $n$-dimensional complete Riemannian manifold and $\gamma$ be a minimal geodesic joining two points of $M$. Assume that
\[ Ric(\gamma',\gamma') \geq (n-1)H + \frac{d\phi}{dt} \]
holds along $\gamma$, where $H$ is a positive constant and $\phi$ is any smooth function satisfying $|\phi| \leq L$. Then $M$ is compact and its diameter is bounded from above by
\[ \mbox{\rm diam(M)} \leq \frac{\pi}{(n-1)H}(L + \sqrt{L^{2} + (n-1)^{2}H}). \]
\end{theorem}
 
 In this paper, we study Myers-type theorems for complete Riemannian manifolds in the context of the $h$–almost Ricci tensor and generalized quasi–Einstein tensor. The following theorem generalizes original Myers theorem via $h$–almost Ricci tensor:

\begin{theorem} \label{thm1}
 Let $M$ be an $n$-dimensional complete Riemannian manifold. Suppose that there exists some positive constant $H>0$ such that the $h$–almost Ricci tensor satisfies
\be \label{eq1}
Ric^{\, h}_{\, V}(\gamma',\gamma') \geq (n-1)H,
\ee
$|V| \leq k_{1}$, $|h| \leq k_{2}$, and $|h'| \leq k_{3}$ for some constants $k_{1}k_{3} < (n-1)H$, $k_{2} > 0$. Then $M$ is compact and
\[\mbox{\rm diam(M)} \leq \frac{1}{((n-1)H - k_{1}k_{3})}\left(2k_{1}k_{2} + \sqrt{4k_{1}^{2}k_{2}^{2} + ((n-1)H - k_{1}k_{3})((n-1)\pi^{2})} \right).\]
\end{theorem}

 As to the theorem \ref{thm1.2} via $h$–almost Ricci tensor, we prove the following compactness theorem:

\begin{theorem} \label{thm2}
 Let $M$ be an $n$-dimensional complete Riemannian manifold. Suppose that there exist some constants $H>0$ and $L \geq 0$ such that for every pair of points in $M$ and minimal geodesic $\gamma$ joining those points, the $h$–almost Ricci tensor satisfies
\be \label{eq2}
Ric^{\, h}_{\, V}(\gamma',\gamma') \geq (n-1)H + \frac{d\phi}{dt},
\ee
where $\phi$ is some smooth function of the arc length satisfying $|\phi| \leq L$ along $\gamma$. If the vector field $V$ and a smooth fuction $h$ satisfy $|V| \leq k_{1}$, $|h| \leq k_{2}$, and $|h'| \leq k_{3}$ for some constants $k_{1}k_{3} < (n-1)H$, $k_{2} > 0$, then $M$ is compact and
\[\mbox{\rm diam(M)} \leq \frac{2(L + k_{1}k_{2})}{((n-1)H - k_{1}k_{3})} + \frac{\sqrt{4(L + k_{1}k_{2})^{2} + ((n-1)H - k_{1}k_{3})((n-1)\pi^{2})}}{((n-1)H - k_{1}k_{3})}. \]
\end{theorem}

\begin{remark}
\emph{By taking $L=0$, theorem \ref{thm2} is reduced to the Myers theorem via $h$–almost Ricci tensor (theorem \ref{thm1}) with the diameter estimate}
\bea
\mbox{\rm diam(M)} \leq \frac{1}{((n-1)H - k_{1}k_{3})}\left(2k_{1}k_{2} + \sqrt{4k_{1}^{2}k_{2}^{2} + ((n-1)H - k_{1}k_{3})((n-1)\pi^{2})} \right).
\eea
\end{remark}

\par\vspace{.5cm}

 Moreover, we prove the Myers theorem of a complete Riemannian manifold with a lower bound on the $h$–almost Ricci tensor under the condition that the smooth function $h$ in the $h$–almost Ricci tensor has at most linear growth in the distance function.

\begin{theorem} \label{thm3}
 Let $M$ be an $n$-dimensional complete Riemannian manifold. Suppose that there exists some positive constant $H>0$ such that the $h$–almost Ricci tensor satisfies
\be \label{eq3}
Ric^{\, h}_{\, V}(\gamma',\gamma') \geq (n-1)H,
\ee
$|V| \leq k_{1}$, $|h| \leq k_{2}(d(x,p) + 1)$, and $|h'| \leq k_{2}$ for some constant $k_{1}k_{2} < \frac{(n-1)H}{2}$, where $d(x,p)$ is the distance function from some fixed $p$ to $x$. Then $M$ is compact.
\end{theorem}

 Furthermore, we prove the Myers-type theorem for a complete Riemannian manifold when $h$ has at most linear growth.

\begin{theorem} \label{thm4}
 Let $M$ be an $n$-dimensional complete Riemannian manifold. Suppose that there exist some constants $H>0$ and $L \geq 0$ such that for every pair of points in $M$ and minimal geodesic $\gamma$ joining those points, the $h$–almost Ricci tensor satisfies
\be \label{eq4}
Ric^{\, h}_{\, V}(\gamma',\gamma') \geq (n-1)H + \frac{d\phi}{dt},
\ee
where $\phi$ is some smooth function of the arc length satisfying $\phi \geq -L$ along $\gamma$. If the vector field $V$ and a smooth function $h$ satisfy $|V| \leq k_{1}$, $|h| \leq k_{2}(d(x,p) + 1)$, and $|h'| \leq k_{2}$ for some constant $k_{1}k_{2} < \frac{(n-1)H}{2}$, where $d(x,p)$ is the distance function from some fixed $p$ to $x$, then $M$ is compact.
\end{theorem}

\par\vspace{.5cm}

 Now we consider generalized quasi–Einstein tensor. First, we study original Myers theorem via generalized quasi–Einstein tensor.

\begin{theorem} \label{thm5}
 Let $M$ be an $n$-dimensional complete Riemannian manifold. Suppose that there exists some positive constant $H >0$ such that a generalized quasi–Einstein tensor satisfies
\be \label{eq5}
Ric_{f}^{\, \mu}(\gamma',\gamma') \geq (n-1)H,
\ee
where $\mu \geq \frac{1}{k_{4}}$ for some positive constant $k_{4}$. Then $M$ is compact and
\[ \mbox{\rm diam(M)} \leq \frac{\pi}{(n-1)H}\sqrt{(n+k_{4}-1)(n-1)H}. \]
\end{theorem}

 Secondly, we prove the following compactness theorem for generalized quasi–Einstein tensor which generalizes theorem \ref{thm1.2}:

\begin{theorem} \label{thm6}
Let $M$ be an $n$-dimensional complete Riemannian manifold. Suppose that there exist some constants $H>0$ and $L \geq 0$ such that for every pair of points in $M$ and minimal geodesic $\gamma$ joining those points, a generalized quasi–Einstein tensor satisfies
\be \label{eq6}
Ric_{f}^{\, \mu}(\gamma',\gamma') \geq (n-1)H + \frac{d\phi}{dt},
\ee
where $\phi$ is some smooth function of the arc length satisfying $|\phi| \leq L$ along $\gamma$. If a smooth function $\mu$ satisfies $\mu \geq \frac{1}{k_{4}}$ for some positive constant $k_{4}$, then $M$ is compact and
\[\mbox{\rm diam(M)} \leq \frac{1}{(n-1)H}\left(2L + \sqrt{4L^{2} + (n-1)H(n+k_{4}-1)\pi^{2}} \right). \]
\end{theorem}

\par\vspace{.5cm}

 Finally, we generalize Theorem \ref{thm1.1} to generalized quasi–Einstein tensors.

\begin{theorem} \label{thm7}
 Let $M$ be an $n$-dimensional complete Riemannian manifold. Suppose that there exists some point $p \in M$ for which every geodesic $\gamma: [0, \infty) \rightarrow M$ emanating from $p$ satisfies
\be \label{eq7}
\int^{\infty}_{0} Ric_{f}^{\, \mu}(\gamma'(s),\gamma'(s))ds = \infty,
\ee
where $\mu \geq \frac{1}{k_{4}}$ for some positive constant $k_{4}$. Then $M$ is compact.
\end{theorem}

 This paper is organized as follows: In Section 2 and Section 3, we prove Myers-type theorems for $h$–almost Ricci tensor. In Section 4, we study Myers-type theorems for generalized quasi–Einstein tensor.

\section{When $h$ is bounded to a constant}
 In this section, we prove Theorem \ref{thm1} and Theorem \ref{thm2}. First, we prove Theorem \ref{thm1}. Let $p,q \in M$ and $\gamma$ be a minimizing unit speed geodesic segment from $p$ to $q$ of length $\ell$. Consider a parallel orthonormal frame $\{E_{1}=E_{1}, E_{2}, \cdots, E_{n}=\gamma' \}$ along $\gamma$ and a smooth function $b \in C^{\infty}([0,\ell])$ such that $b(0) = b(\ell) = 0$. By the index form, we have
\[ \sum_{i=1}^{n-1} I(bE_{i},bE_{i}) = \int^{\ell}_{0} (n-1)(b')^{2} - b^{2}Ric(\gamma',\gamma') dt, \]
where $I$ denotes the index form of $\gamma$. Using the assumption (\ref{eq1}) and definition of $h$-almost Ricci tensor, we get
\[\sum_{i=1}^{n-1} I(bE_{i},bE_{i}) \leq (n-1)\int^{\ell}_{0} (b')^{2} dt - (n-1)H\int^{\ell}_{0} b^{2} dt + \int^{\ell}_{0} \frac{b^{2}h}{2}\mathcal{L}_{V}g(\gamma',\gamma') dt.\]
Note that
\bea
\int^{\ell}_{0} \frac{b^{2}h}{2}\mathcal{L}_{V}g(\gamma',\gamma') dt &=& \int^{\ell}_{0} b^{2}h\frac{d}{dt}\langle V,\gamma' \rangle dt \\
&=& \int^{\ell}_{0} \frac{d}{dt}(b^{2}h\langle V,\gamma' \rangle)dt - 2\int^{\ell}_{0} bb'h\langle V,\gamma' \rangle dt - \int^{\ell}_{0}b^{2}h'\langle V,\gamma' \rangle dt \\
&\leq& 2\int^{\ell}_{0} |bb'||h|\left|\langle V,\gamma' \rangle \right| dt + \int^{\ell}_{0}|b^{2}||h'|\left|\langle V,\gamma' \rangle \right| dt \\
&\leq& 2k_{1}k_{2}\int^{\ell}_{0} |bb'| dt + k_{1}k_{3}\int^{\ell}_{0} |b^{2}| dt.
\eea
So we have
\bea
\sum_{i=1}^{n-1} I(bE_{i},bE_{i}) &\leq& (n-1)\int^{\ell}_{0} (b')^{2} dt - (n-1)H\int^{\ell}_{0} b^{2} dt + 2k_{1}k_{2}\int^{\ell}_{0} |bb'| dt \\
&& + k_{1}k_{3}\int^{\ell}_{0} |b^{2}| dt.
\eea
By setting the function $b$ to be $b(t) = sin(\frac{\pi t}{\ell})$, we have $b'(t) = \frac{\pi}{\ell}cos(\frac{\pi t}{\ell})$. Thus,
\bea
\sum_{i=1}^{n-1} I(bE_{i},bE_{i}) &\leq& \frac{(n-1)\pi^{2}}{\ell^{2}} \int^{\ell}_{0} cos^{2}(\frac{\pi t}{\ell}) dt - (n-1)H\int^{\ell}_{0}sin^{2}(\frac{\pi t}{\ell})dt \\
&&+ \frac{k_{1}k_{2} \pi}{\ell}\int^{\ell}_{0} \left|sin(\frac{2\pi t}{\ell})\right| dt + k_{1}k_{3}\int^{\ell}_{0}\left|sin^{2}(\frac{\pi t}{\ell})\right|dt \\
&\leq& \frac{(n-1)\pi^{2}}{2\ell} - \frac{(n-1)H\ell}{2} + 2k_{1}k_{2} + \frac{k_{1}k_{3}\ell}{2} \\
&=& -\frac{1}{2\ell}\left( ((n-1)H - k_{1}k_{3})\ell^{2} - 4k_{1}k_{2}\ell - (n-1)\pi^{2}\right ).
\eea
Since $\gamma$ is a minimizing geodesic, we must take
\[ ((n-1)H - k_{1}k_{3})\ell^{2} - 4k_{1}k_{2}\ell - (n-1)\pi^{2} \leq 0. \]
If $k_{1}k_{3} < (n-1)H$, then $(n-1)H - k_{1}k_{3} > 0$, so we have
\[ \ell \leq \frac{1}{(n-1)H - k_{1}k_{3}} \left(2k_{1}k_{2} + \sqrt{4k_{1}^{2}k_{2}^{2} + ((n-1)H - k_{1}k_{3})((n-1)\pi^{2})} \right ). \]
Hence, $M$ is compact and
\[ \mbox{\rm diam(M)} \leq \frac{1}{(n-1)H - k_{1}k_{3}} \left(2k_{1}k_{2} + \sqrt{4k_{1}^{2}k_{2}^{2} + ((n-1)H - k_{1}k_{3})((n-1)\pi^{2})} \right ). \]
This completes the proof.
\hfill\fbox{}\par\vspace{.5cm}

Second, we prove Theorem \ref{thm2}. Its proof is similar to the previous proof; thus, the setting is identical. The index form implies
\[ \sum_{i=1}^{n-1} I(bE_{i},bE_{i}) = \int^{\ell}_{0} (n-1)(b')^{2} - b^{2}Ric(\gamma',\gamma') dt. \]
From (\ref{eq2}), we have
\bea
\sum_{i=1}^{n-1} I(bE_{i},bE_{i}) &\leq& (n-1)\int^{\ell}_{0}(b')^{2} dt - (n-1)H\int^{\ell}_{0}b^{2} dt -\int^{\ell}_{0}b^{2}\frac{d\phi}{dt} dt \\
&& + \int^{\ell}_{0} \frac{1}{2}b^{2}h\mathcal{L}_{V}g(\gamma',\gamma') dt.
\eea
Note that
\[ \int^{\ell}_{0} \frac{1}{2}b^{2}h\mathcal{L}_{V}g(\gamma',\gamma') dt \leq 2k_{1}k_{2}\int^{\ell}_{0} |bb'| dt + k_{1}k_{3}\int^{\ell}_{0} |b^{2}| dt \]
and
\bea
-\int^{\ell}_{0} b^{2}\frac{d\phi}{dt} dt &=& -\left(\int^{\ell}_{0} \frac{d}{dt}(b^{2}\phi) dt - 2\int^{\ell}_{0} bb'\phi \, dt \right) \\
&\leq& 2\int^{\ell}_{0}|bb'||\phi| dt 
\leq 2L\int^{\ell}_{0}|bb'| dt.
\eea 
Thus,
\bea
\sum_{i=1}^{n-1} I(bE_{i},bE_{i}) &\leq& (n-1)\int^{\ell}_{0} (b')^{2} dt - (n-1)H\int^{\ell}_{0} b^{2} dt + 2L\int^{\ell}_{0}|bb'|dt  \\
&& + 2k_{1}k_{2}\int^{\ell}_{0}|bb'|dt + k_{1}k_{3}\int^{\ell}_{0}|b^{2}| dt \\
&=& (n-1)\int^{\ell}_{0} (b')^{2} dt - (n-1)H\int^{\ell}_{0} b^{2} dt + 2(L + k_{1}k_{2})\int^{\ell}_{0}|bb'|dt \\
&& + k_{1}k_{3}\int^{\ell}_{0}|b^{2}| dt.
\eea
If the fuction $b$ is taken to be $b(t) = sin(\frac{\pi t}{\ell})$, then we obtain
\bea
\sum_{i=1}^{n-1} I(bE_{i},bE_{i}) &\leq& \frac{(n-1)\pi^{2}}{\ell^{2}} \int^{\ell}_{0} cos^{2}(\frac{\pi t}{\ell}) \, dt - (n-1)H\int^{\ell}_{0}sin^{2}(\frac{\pi t}{\ell}) \, dt \\
&& + \frac{\pi (L + k_{1}k_{2})}{\ell}\int^{\ell}_{0} \left|sin(\frac{2\pi t}{\ell}) \right| dt + k_{1}k_{3}\int^{\ell}_{0}\left|sin^{2}(\frac{\pi t}{\ell}) \right|dt \\
&=& \frac{(n-1)\pi^{2}}{2\ell} - \frac{(n-1)H\ell}{2} + 2(L + k_{1}k_{2}) + \frac{k_{1}k_{3}\ell}{2} \\
&=& -\frac{1}{2\ell}\left(((n-1)H - k_{1}k_{3})\ell^{2} - 4\ell(L + k_{1}k_{2}) - (n-1)\pi^{2}\right).
\eea
Since $\gamma$ is a minimizing geodesic, we must take
\[ ((n-1)H - k_{1}k_{3})\ell^{2} - 4\ell(L + k_{1}k_{2}) - (n-1)\pi^{2} \leq 0. \]
If $k_{1}k_{3} < (n-1)H$, then $(n-1)H - k_{1}k_{3} > 0$, so we get
\[ \ell \leq \frac{2(L + k_{1}k_{2})}{((n-1)H - k_{1}k_{3})} + \frac{\sqrt{4(L + k_{1}k_{2})^{2} + ((n-1)H - k_{1}k_{3})((n-1)\pi^{2})}}{((n-1)H - k_{1}k_{3})}. \]
This proves Theorem \ref{thm2}.
\hfill\fbox{}\par\vspace{.5cm}

\section{When $h$ has at most linear growth}
 In this section, we prove Myers-type theorems (Theorem \ref{thm3} and \ref{thm4}) under the condition that $h$ has at most linear growth. Our proofs use the Riccati inequality. 
\\
 
 First, we prove Theorem \ref{thm3}. Suppose that $M$ is non-compact. Fix a point $p \in M$, there exists a unit speed ray $\gamma(t)$ starting from $p$ satisfying $\gamma(0) = p$. For every $t>0$, the mean curvature function $m(t)$ defined by $m(x)= \Delta r(x)$, where $r(x) = d(p,x)$ satisfies the Riccati inequality:
\bea
Ric(\gamma',\gamma') \leq -m'(t) - \frac{1}{n-1}m^{2}(t).
\eea
By adding $\frac{h(t)}{2}\mathcal{L}_{V}g(\gamma',\gamma')$ to both sides of this inequality, we have
\bea
Ric(\gamma'(t),\gamma'(t)) + \frac{h(t)}{2}\mathcal{L}_{V}g(\gamma'(t),\gamma'(t)) &\leq& -m'(t) - \frac{1}{n-1}m^{2}(t) \\ 
&& + \frac{h(t)}{2}\mathcal{L}_{V}g(\gamma'(t),\gamma'(t)) \\
\eea
From (\ref{eq3}) and definition of $Ric_{\, V}^{\, h}$, we obtain
\be \label{eq3.1}
(n-1)H \leq -m'(t) - \frac{1}{n-1}m^{2}(t) + \frac{h(t)}{2}\mathcal{L}_{V}g(\gamma'(t),\gamma'(t)).
\ee
Integrating both sides of (\ref{eq3.1}), we have
\bea
\int^{t}_{1}(n-1)H \, ds &\leq& -m(t) + m(1) -\frac{1}{n-1}\int^{t}_{1}m^{2}(s) \, ds \\
&& + \int^{t}_{1} \frac{h(s)}{2}\mathcal{L}_{V}g(\gamma'(s),\gamma'(s)) \, ds.
\eea
Note that
\bea
\int^{t}_{1} \frac{h(s)}{2}\mathcal{L}_{V}g(\gamma'(s),\gamma'(s)) \, ds &=& \int^{t}_{1}h(s)\frac{d}{ds}\langle V,\gamma' \rangle(s) \, ds \\
&=& \int^{t}_{1}\frac{d}{ds}\left(h(s)\langle V,\gamma' \rangle(s)\right) ds - \int^{t}_{1}h'(s)\langle V,\gamma' \rangle(s) \, ds \\
&=& h(t)\langle V,\gamma' \rangle(t) - h(1)\langle V,\gamma' \rangle(1) - \int^{t}_{1}h'(s)\langle V,\gamma' \rangle(s) \, ds.
\eea
So we get
\bea
(n-1)Ht - (n-1)H &\leq& -m(t) + m(1) - \frac{1}{n-1}\int^{t}_{1}m^{2}(s) \, ds + h(t)\langle V,\gamma' \rangle(t) \\
&& - h(1)\langle V,\gamma' \rangle(1) - \int^{t}_{1}h'(s)\langle V,\gamma' \rangle(s) \, ds \\
&\leq& -m(t) + m(1) - \frac{1}{n-1}\int^{t}_{1}m^{2}(s) \, ds + |h(t)||\langle V,\gamma' \rangle(t)| \\
&& - h(1)\langle V,\gamma' \rangle(1) + \int^{t}_{1}|h'(s)||\langle V,\gamma' \rangle(s)| \, ds \\
&\leq& -m(t) + m(1) - \frac{1}{n-1}\int^{t}_{1}m^{2}(s) \, ds + 2k_{1}k_{2}t - h(1)\langle V,\gamma' \rangle(1).
\eea
It follows that
\bea
-m(t) - \frac{1}{n-1}\int^{t}_{1}m^{2}(s) \, ds &\geq& ((n-1)H - 2k_{1}k_{2})t + C_{1},
\eea
where $C_{1}:= + h(1)\langle V,\gamma' \rangle(1) - (n-1)H - m(1)$.
Since $k_{1}k_{2} < \frac{(n-1)H}{2}$, there exists $t_{1}>1$ such that for all $t \geq t_{1}$, we have
\be \label{eq3.2}
-m(t) - \frac{1}{n-1}\int^{t}_{1}m^{2}(s) \, ds \geq 2.
\ee
Now we consider the increasing sequence $\{t_{\ell}\}$ defined by 
\[ t_{\ell + 1} = t_{\ell} + (n-1)2^{1 - \ell}, \,\, {\rm for }\ \ell \geq 1. \]
Note that $\{t_{\ell}\}$ converges to $T:= t_{1} + 2(n-1)$ as $\ell \rightarrow \infty$.

 We claim that $-m(t) \geq 2^{\ell}$ for all $t \geq t_{\ell}$. To prove this claim, we use induction. If $\ell = 1$, the claim is trivially true from the inequality in (\ref{eq3.2}). Now, for all $t \geq t_{\ell + 1}$, we have
\bea
-m(t) &\geq& 2 + \frac{1}{n-1}\int^{t}_{1} m^{2}(s) \, ds \\
&\geq& \frac{1}{n-1}\int^{t_{\ell +1}}_{t_{\ell}} m^{2}(s) ds \\
&\geq& \frac{1}{n-1}2^{2\ell}(t_{\ell +1} - t_{\ell}) = 2^{\ell + 1}.
\eea
Hence, the claim is true for all $t \geq t_{\ell +1}$. Therefore,
\[ \lim_{\ell \rightarrow \infty}-m(t_{\ell}) = -m(T) \geq \lim_{\ell \rightarrow \infty}2^{\ell + 1}. \]
This contradicts the smoothness of $m(t)$, which completes the proof of Theorem \ref{thm3}.
\hfill\fbox{}\par\vspace{.5cm}

 Second, we prove Theorem \ref{thm4}. Setting is the same as the above proof. We have
\bea
Ric(\gamma'(t),\gamma'(t)) + \frac{h(t)}{2}\mathcal{L}_{V}g(\gamma'(t),\gamma'(t)) &\leq& -m'(t) - \frac{1}{n-1}m^{2}(t) \\ 
&& + \frac{h(t)}{2}\mathcal{L}_{V}g(\gamma'(t),\gamma'(t)).
\eea
From (\ref{eq4}), we obtain
\be \label{eq3.3}
(n-1)H + \frac{d\phi}{dt} \leq -m'(t) - \frac{1}{n-1}m^{2}(t) + \frac{h(t)}{2}\mathcal{L}_{V}g(\gamma'(t),\gamma'(t)).
\ee
Integrating both sides of (\ref{eq3.3}) from $1$ to $t$, we get
\bea
(n-1)Ht - (n-1)H + \phi(t) - \phi(1) &\leq& -m(t) + m(1) - \frac{1}{n-1}\int^{t}_{1}m^{2}(s) \, ds \\
&& + \int^{t}_{1} \frac{h(s)}{2}\mathcal{L}_{V}g(\gamma'(s),\gamma'(s)) \, ds.
\eea
Note that
\bea
\int^{t}_{1} \frac{h(s)}{2}\mathcal{L}_{V}g(\gamma'(s),\gamma'(s)) \, ds &=& \int^{t}_{1}h(s)\frac{d}{ds}\langle V,\gamma' \rangle(s) \, ds \\
&=& h(t)\langle V,\gamma' \rangle(t) - h(1)\langle V,\gamma' \rangle(1)  - \int^{t}_{1}h'(s)\langle V,\gamma' \rangle(s) \, ds \\ 
&\leq& + |h(t)||\langle V,\gamma' \rangle(t)| - h(1)\langle V,\gamma' \rangle(1) + \int^{t}_{1}|h'(s)||\langle V,\gamma' \rangle(s)| \, ds \\
&\leq& 2k_{1}k_{2}t - h(1)\langle V,\gamma' \rangle(1).
\eea
If $\phi \geq -L$, then we have
\bea
(n-1)Ht - (n-1)H - L - \phi(1) &\leq& -m(t) + m(1) - \frac{1}{n-1}\int^{t}_{1}m^{2}(s) \, ds \\
&& + 2k_{1}k_{2}t - h(1)\langle V,\gamma' \rangle(1).
\eea
Thus,
\bea
-m(t) - \frac{1}{n-1}\int^{t}_{1}m^{2}(s) \, ds &\geq& (n-1)Ht - (n-1)H - L - \phi(1) \\
&& -m(1) - 2k_{1}k_{2}t + h(1)\langle V,\gamma' \rangle(1).
\eea
Let $C_{2}:= - (n-1)H - L - \phi(1) + h(1)\langle V,\gamma' \rangle(1)$. Then we obtain
\bea
-m(t) - \frac{1}{n-1}\int^{t}_{1}m^{2}(s) \, ds \geq t\left((n-1)H -2k_{1}k_{2}t \right) + C_{2}.
\eea
Since $k_{1}k_{2} < \frac{(n-1)H}{2}$, the above inequality implies that there exists $t_{1} > 1$ such that for all $t \geq t_{1}$, we get
\[ -m(t) - \frac{1}{n-1}\int^{t}_{1}m^{2}(s) \, ds \geq 2. \]
Now we can complete this proof by using the same argement as in the proof of Theorem \ref{thm3}. 

\hfill\fbox{}\par\vspace{.5cm}

\section{generalized quasi–Einstein tensor}
 In this section, we prove Theorem \ref{thm5}, \ref{thm6}, and \ref{thm7}. First, we prove Theorem \ref{thm5}. The setting is identical with proof of Theorem \ref{thm1}. By the index form, we have
\[ \sum_{i=1}^{n-1} I(bE_{i},bE_{i}) = \int^{\ell}_{0} (n-1)(b')^{2} - b^{2}Ric(\gamma',\gamma') \, dt. \]
By the assumption (\ref{eq5}) and definition of generalized quasi–Einstein tensor (\ref{eq1.6}), we obtain
\bea
\sum_{i=1}^{n-1} I(bE_{i},bE_{i}) &\leq& (n-1)\int^{\ell}_{0}(b')^{2} \, dt - (n-1)H\int^{\ell}_{0}b^{2} \, dt + \int^{\ell}_{0} b^{2}\mbox{\rm Hess}f(\gamma',\gamma') \, dt \\
&&- \int^{\ell}_{0}b^{2}\mu \, df \otimes df(\gamma',\gamma') \, dt.
\eea
Note that
\bea
\int^{\ell}_{0} b^{2}\mbox{\rm Hess}f(\gamma',\gamma') \, dt &=& \int^{\ell}_{0}b^{2}\frac{d}{dt}\langle\nabla f, \gamma' \rangle(t) \, dt \\
&=& \int^{\ell}_{0}\frac{d}{dt}\left(b^{2}\langle\nabla f, \gamma' \rangle(t) \right)dt - 2\int^{\ell}_{0}bb'\langle\nabla f, \gamma' \rangle(t) \, dt \\
&\leq& 2\int^{\ell}_{0}|b'||b\langle\nabla f, \gamma' \rangle(t)| \, dt \\
&\leq& 2 \left(\int^{\ell}_{0}\frac{1}{\mu}(b')^{2}dt \right)^{\frac{1}{2}}\left(\int^{\ell}_{0}\mu\left(b \, \langle\nabla f, \gamma' \rangle(t)\right)^{2}dt \right)^{\frac{1}{2}} \\
&\leq& \int^{\ell}_{0} \frac{1}{\mu}(b')^{2}dt + \int^{\ell}_{0}\mu\left(b \, \langle\nabla f, \gamma' \rangle(t)\right)^{2}dt.
\eea
So we have
\bea
\sum_{i=1}^{n-1} I(bE_{i},bE_{i}) &\leq& (n-1)\int^{\ell}_{0}(b')^{2} \, dt - (n-1)H\int^{\ell}_{0}b^{2} \, dt + \int^{\ell}_{0} \frac{1}{\mu}(b')^{2}dt \\
&\leq& (n-1)\int^{\ell}_{0}(b')^{2} \, dt - (n-1)H\int^{\ell}_{0}b^{2} \, dt + k_{4}\int^{\ell}_{0}(b')^{2}dt.
\eea
If $b(t) = sin(\frac{\pi t}{\ell})$, then we obtain
\bea
\sum_{i=1}^{n-1} I(bE_{i},bE_{i}) &\leq& \frac{\pi^{2}(n-1)}{2\ell} - \frac{(n-1)H\ell}{2} + \frac{\pi^{2}k_{4}}{2\ell} \\
&\leq& -\frac{1}{2\ell}\left((n-1)H\ell^{2} -\pi^{2}(n + k_{4} - 1)\right).
\eea
We must take
\[ (n-1)H\ell^{2} -\pi^{2}(n + k_{4} - 1) \leq 0. \]
Hence, we get
\[ \ell \leq \frac{\pi\sqrt{(n-1)H(n+k_{4}-1)}}{(n-1)H}. \]
So we complete the proof of Theorem \ref{thm5}.
\hfill\fbox{}\par\vspace{.5cm}

 Now we prove Theorem \ref{thm6}. Setting is the same as the above proof. From the index form, generalized quasi–Einstein tensor (\ref{eq1.6}), and the assumption (\ref{eq6}), we have
\bea
\sum_{i=1}^{n-1} I(bE_{i},bE_{i}) &\leq& (n-1)\int^{\ell}_{0}(b')^{2} \, dt - (n-1)H\int^{\ell}_{0}b^{2} \, dt - \int^{\ell}_{0}b^{2}\frac{d\phi}{dt} dt  \\
&&+ \int^{\ell}_{0} b^{2}\mbox{\rm Hess}f(\gamma',\gamma') \, dt - \int^{\ell}_{0}b^{2}\mu \, df \otimes df(\gamma',\gamma') \, dt.
\eea
Note that
\bea
\int^{\ell}_{0} b^{2}\mbox{\rm Hess}f(\gamma',\gamma') \, dt \leq \int^{\ell}_{0} \frac{1}{\mu}(b')^{2}dt + \int^{\ell}_{0}\mu\left(b \, \langle\nabla f, \gamma' \rangle(t)\right)^{2}dt,
\eea
and
\bea
-\int^{\ell}_{0} b^{2}\frac{d\phi}{dt} dt \leq 2L\int^{\ell}_{0}|bb'| dt. 
\eea
It follows that
\bea
\sum_{i=1}^{n-1} I(bE_{i},bE_{i}) &\leq& (n-1)\int^{\ell}_{0}(b')^{2} \, dt - (n-1)H\int^{\ell}_{0}b^{2} \, dt + 2L\int^{\ell}_{0}|bb'| dt + \int^{\ell}_{0} \frac{1}{\mu}(b')^{2}dt \\
&\leq& (n-1)\int^{\ell}_{0}(b')^{2} \, dt - (n-1)H\int^{\ell}_{0}b^{2} \, dt + 2L\int^{\ell}_{0}|bb'| dt + k_{4}\int^{\ell}_{0}(b')^{2}dt.
\eea
If $b(t) = sin(\frac{\pi t}{\ell})$, then we have
\[ \sum_{i=1}^{n-1} I(bE_{i},bE_{i}) \leq -\frac{1}{2\ell}\left((n-1)H\ell^{2} - 4L\ell - \pi^{2}(n + k_{4} -1)\right). \]
So we must take
\[ (n-1)H\ell^{2} - 4L\ell - \pi^{2}(n + k_{4} -1) \leq 0. \]
Therefore, we obtain
\[ \ell \leq \frac{1}{(n-1)H}\left(2L + \sqrt{4L^{2} + \pi^{2}(n+k_{4}-1)(n-1)H} \right). \]
This proves Theorem \ref{thm6}.
\hfill\fbox{}\par\vspace{.5cm}

 Finally, we prove Theorem \ref{thm7}. Before the proof, we fix several notation. Let $(M,g,e^{-f}dv_{g})$ be a smooth metric measure space on an $n$-dimensional complete Riemannian manifold $M$. For the measure $e^{-f}dv_{g}$, the $f$-mean curvature is $m_{f} = m - \partial_{r}f$, where $m$ is the mean curvature of the geodesic sphere with inward pointing normal vector. Then the $f$-Laplacian is defined by $\Delta_{f} := \Delta - \langle\nabla f, \nabla \rangle$. Note that $m_{f} = \Delta_{f}(r)$ and $m = \Delta(r)$, where $r$ is the distance function.

 Fix a point $p \in M$ and take a unit speed ray $\gamma = \gamma(s)$ emanating from $p$ satisfying $\gamma(0) = p$. Let $r(x) = d(x,p)$ be the distance between $x$ and $p$. By Bochner formula and Schwarz inequality, the distance function $r$ satisfies the Riccati inequality
\[ \frac{\partial}{\partial r}(\Delta r) \leq -\frac{1}{n-1}(\Delta r)^{2} - Ric(\nabla r,\nabla r) .\]
We know that $\Delta_{f}(r) = \Delta r - \langle\nabla f, \nabla r \rangle$. So we have
\[ \frac{\partial}{\partial r}(\Delta_{f}(r)) \leq -\frac{1}{n-1}\left(\Delta_{f}(r) + \langle\nabla f, \nabla r \rangle\right)^{2} - Ric(\nabla r, \nabla r) - \mbox{\rm Hess}f(\nabla r, \nabla r). \]
Recall the elementary inequality
\be \label{eq4.1}
(a+b)^{2} \geq \frac{1}{\alpha+1}a^{2} - \frac{1}{\alpha}b^{2}
\ee
for $\alpha > 0$. By (\ref{eq4.1}), we obtain
\[ \frac{\partial}{\partial r}(\Delta_{f}(r)) \leq -\frac{(\Delta_{f}(r))^{2}}{(n-1)(\alpha + 1)} + \frac{\langle\nabla f, \nabla r \rangle^{2}}{(n-1)\alpha} - Ric(\nabla r, \nabla r) - \mbox{\rm Hess}f(\nabla r, \nabla r). \]
Let $(n-1)\alpha = k_{4}$. Then we have
\bea
\frac{\partial}{\partial r}(\Delta_{f}(r)) &\leq&  -\frac{(\Delta_{f}(r))^{2}}{k_{4}+ n-1}+ \frac{\langle\nabla f, \nabla r \rangle^{2}}{k_{4}} - Ric(\nabla r, \nabla r) - \mbox{\rm Hess}f(\nabla r, \nabla r) \\
&\leq& -\frac{(\Delta_{f}(r))^{2}}{k_{4}+ n-1}+ \mu\langle\nabla f, \nabla r \rangle^{2} - Ric(\nabla r, \nabla r) - \mbox{\rm Hess}f(\nabla r, \nabla r) \\ 
&=& -\frac{(\Delta_{f}(r))^{2}}{k_{4}+ n-1} - Ric_{f}^{\, \mu}(\nabla r, \nabla r).
\eea
It follows that
\be \label{eq4.2}
Ric_{f}^{\, \mu}(\gamma',\gamma') \leq -m'_{f}(s) - \frac{m_{f}^{2}(s)}{k_{4}+n-1},
\ee
where $m_{f}(s) = (\Delta_{f}\, r)(\gamma(s))$. \\
Integrating both sides of (\ref{eq4.2}), we have
\[ \lim_{t \rightarrow \infty}\int^{t}_{1} Ric_{f}^{\, \mu}(\gamma'(s),\gamma'(s))ds \leq \lim_{t \rightarrow \infty}\int^{t}_{1}\left(-m'_{f}(s) - \frac{m_{f}^{2}(s)}{k_{4}+n-1} \right) ds. \]
If \[\lim_{t \rightarrow \infty}\int^{t}_{0} Ric_{f}^{\, \mu}(\gamma'(s),\gamma'(s))ds = \infty, \] then we obtain
\[ \lim_{t \rightarrow \infty}\left(-m_{f}(t) - \frac{1}{k_{4}+n-1}\int^{t}_{1}m_{f}^{2}(s) ds \right) = \infty.\]
So there exists $t_{1} > 1$ such that for all $t \geq t_{1}$, we get
\[ -m_{f}(t) - \frac{1}{k_{4}+n-1}\int^{t}_{1}m_{f}^{2}(s) ds > 2. \]
Using the same argument as in the proof of Theorem \ref{thm3}, we may complete the proof.
\hfill\fbox{}\par
\end{document}